\documentclass [12pt] {article}
\usepackage{amsfonts}
\usepackage{amssymb}

\newcommand{\qed} {\hspace {0.1in} \rule {1.5mm} {3.5mm}}

\newtheorem{lemma}{Lemma}[section]
\newtheorem{corollary}[lemma]{Corollary}
\newtheorem{theorem}{Theorem}

\newtheorem{conjecture}{Conjecture}
\newtheorem{proposition}[lemma]{Proposition}
\newtheorem{definition}[lemma]{Definition}

\def\ai{\alpha\in I}
\def\un{\underline}
\def\a{_\alpha}

\def\map{\mathbb Map}

\def\ke{\mbox{Ker}\,}
\def\ra{\mbox{Ran}\,}

\def\r{\mbox{rank}_D\,}

\def\dim{{\rm dim}}

\def\<{\langle}
\def\>{\rangle}
\def\proof{\smallskip\noindent{\bf Proof:} }

\def\bR{{\mathbb R}}
\def\bN{{\mathbb N}}
\def\bC{{\mathbb C}}

\def\en{\mbox{End}}

\def\to{\rightarrow}

\title{Sofic groups and direct finiteness}
\author{{\sc G\'abor Elek and Endre Szab\'o}
\cr Mathematical Institute of
the Hungarian Academy of Sciences\cr P.O. Box 127, H-1364 Budapest, Hungary\cr
elek@renyi.hu}
\date{}
\begin{document}

\maketitle
\noindent{\bf Abstract.} We construct an analogue of von Neumann's affiliated
algebras for sofic group algebras over arbitrary fields. Consequently, we
settle Kaplansky's direct finiteness conjecture for sofic groups. 
\vskip 0.2in
\newpage
\section{Introduction}
The following conjecture is due to Kaplansky.
\begin{conjecture}\label{conj1}
For any
group $G$ and commutative field $K$, the group algebra $K(G)$ is
directly finite. That is $ab=1$ in $K(G)$ implies $ba=1$.
\end{conjecture}

Recently Ara, O'Meara and Perera \cite{AMP} settled the conjecture
for residually amenable groups even in the case of group
algebras $D(G)$, where $D$ is a division ring. They also proved that
$Mat_{n\times n} (D(G))$ is directly finite as well.

\noindent
It is important to note that Conjecture~\ref{conj1}
holds for any group $G$ if $D=\bC$, the complex field \cite{Kap}.
Indeed,
$\bC(G)$ is a subalgebra of the von Neumann algebra $N(G)$.
The algebra $N(G)$ always satisfies the Ore-condition with respect
to its non-zero divisors. Hence one can consider its classical
ring of fractions $U(G)$. The algebra $U(G)$ is the so-called {\em affiliated
algebra} of $G$ and it is a continuous von Neumann-regular ring \cite{Good},
\cite{Lueck}, hence it is known to be directly finite.

\noindent
Let us turn to another conjecture due to Gottschalk \cite{Got}.
\begin{conjecture}\label{conj2}
Let $G$ be a countable group and $X$  a finite set.
Consider the compact metrizable space $X^G$ of $X$-valued functions
on $G$ equipped with the product topology. Let
$f:X^G\to X^G$ be a continuous map that commutes with the natural
right $G$-action. Then if $f$ is injective, it is surjective as well.
\end{conjecture}

\noindent
In \cite{Gro} Gromov proved Gottschalk's conjecture in the
case of {\em sofic groups} (the name ``sofic groups'' was coined
by Weiss \cite{Wei}). We shall review the definitions and basic
properties of sofic groups in Section~\ref{secsofic}, nevertheless,
let us note that residually amenable groups are sofic groups as well,
and in \cite{ES} we constructed sofic groups that are not residually amenable. 
On the other hand, there seems to be no example yet of a group
which is not sofic.

\noindent
Let us observe that Conjecture~\ref{conj2} implies Conjecture~\ref{conj1}
for finite fields $F$. Indeed, it is enough to prove
Conjecture~\ref{conj1} for countable groups.
Then any element $a$ of the group algebra
$F(G)$ induces a continuous linear map on $F^G$ commuting with the
right $G$-action. Simply, $a$ acts as convolution on the left.
Then $F(G)$ can be identified with the dense
set of elements in $F^G$ having only finitely many non-zero values.
If $ab=1$ on this dense subset then $ab$ must be equal to the identity
on the whole $F^G$. Therefore $b$ is injective and thus it is
a bijective continuous map by our assumption. Consequently, $a$ is
the inverse of $b$, thus $ba=1$.

\noindent

The goal of this paper is to replace the notion of the affiliated 
algebras of  complex group algebras
with something similar for group algebras of sofic groups 
over arbitrary division rings.
\begin{theorem}\label{main-theorem}
Let $G$ be a sofic
group  and
let $D$ be a division ring. Then $D(G)$ can be embedded into an
indecomposable, continuous von Neumann-regular ring $R(G)$.

\end{theorem}

Therefore we extend the result of Ara, O'Meara and Perera to 
the class of sofic groups.


\section{Continuous von Neumann regular rings}
In this section we give a brief summary of the theory of continuous
von Neumann rings, based upon the monograph of Goodearl \cite{Good}.
Recall that a unital ring $R$ is {\em von Neumann regular} if for any
$x\in R$ there exists $y\in R$ such that $xyx=x$. It is equivalent to
say that any finitely generated right ideal of $R$ can be generated by
one single idempotent.
A ring $R$ is called {\em unit-regular} if for any $x\in R$ there exists
a unit $y$ such that $xyx=x$. A ring is called {\em directly finite} if
$xy=1$ implies $yx=1$ and it is called {\em stably finite} if $Mat_{n\times n}(R)$
is directly finite for all $n\geq 1$. Any unit-regular ring is necessarily
stably finite. 
A lattice $L$ is called a {\em continuous geometry} 
if it is modular, complete, complemented and
$$a\wedge(\vee_{\alpha\in I}b_\alpha)=\vee _{\alpha\in I}(a\wedge b_\alpha)$$
$$a\vee(\wedge_{\alpha\in I}b_\alpha)=\wedge _{\alpha\in I}(a\vee b_\alpha)$$
for any linearly ordered subset $\{b_\alpha\}_{\alpha\in I}\subset L$.
A von Neumann regular ring is called {\em continuous} if the lattices of both its
right ideals and left ideals are  continuous geometries.
The continuous von Neumann regular rings are unit-regular, hence they are
stably finite as well. Division rings and matrix rings over division rings
are the simplest examples of continuous von Neumann regular rings. The
first simple continuous von Neumann rings which do not satisfy
either the ascending or the descending chain condition had already been
 constructed by John von Neumann \cite{JVN}. The following proposition
 summarizes
what we need to know about such rings.
\begin{proposition}\label{dimension}
If $R$ is a simple continuous, von Neumann regular ring that
does not satisfy either the ascending or the descending chain condition,
then
\begin{itemize}
\item $K_0(R)=\bR$, in fact there exists a unique non-negative real-valued dimension function
$\dim_R$ on the set of finitely generated projective right modules over $R$
taking all non-negative values such that
\begin{enumerate}
\item $\dim_R(R)=1\,.$
\item $\dim_R(0)=0\,.$
\item $\dim_R(A\oplus B)=\dim_R(A)+\dim_R(B)\,.$
\item $\dim_R(A)=\dim_R(B)$ if and only if $A\simeq B$.
\end{enumerate}
\item If $A,B$ are finitely generated submodules of a projective module
      then so is $A\cap B$, and
      $$\dim_R(A\cap B)+\dim_R(A+B)=\dim_R(A)+\dim_R(B)\,.$$
\item If $a\in R$ then $Ann(a)=\{x\in R:\, ax=0\}$
is a principal right ideal and
$$\dim_R(Ann(a))+\dim_R(aR)=1\,.$$
\item If $A\le B$ are finitely generated projective modules
and $\dim_R(A)=\dim_R(B)$ then $A=B$.
\end{itemize}
\end{proposition}

\section{A pseudo-rank function on the direct product of matrix rings}
\label{omega-stuff}
Let $I$ be a set and let $\{A_\alpha\}_{\alpha\in I}$ be finitely generated
right $D$-modules, where $D$ is a division ring.
Consider the direct product $E=\prod_{\alpha\in I}\,\en_D(A_\alpha)$, where
$\en_D(A_\alpha)$ is the endomorphism ring of $A_\alpha$.
This ring is directly finite, von Neumann regular, right and left
self-injective. Now we recall the notion of a pseudo-rank function
\cite{Good}.
\begin{definition}
A pseudo-rank function on a von Neumann regular ring $R$ is a map
$N:R\to [0,1]$ such that
\begin{description}
\item[(a)] $ N(1)=1,\,N(0)=0.$
\item[(b)] $N(xy)\leq N(x), N(xy)\leq N(y)$, for all $x,y\in \bR$.
\item[(c)] $N(e+f)=N(e)+N(f)$ for all orthogonal idempotents $e,f\in \bR$.
\end{description} \end{definition}
Before stating our proposition let us recall the notion of ultralimit as well.
Let $\omega$ be an ultrafilter on the set $I$. Then $\lim_\omega$ is the
unique real valued functional on the space of bounded real sequences
$\{a_\alpha\}_{\alpha\in I}$ such that:

\noindent
If $lim_\omega (\{a_\alpha\}_{\alpha\in I})=t$, then
for any $\epsilon>0$:
$$\{\alpha\in I\,|\, a_\alpha\in [t-\epsilon, t+\epsilon]\}\in \omega\,.$$
Note that $\lim_\omega(\{a_\alpha+b_\alpha\}_{\alpha\in I})=\lim_\omega
(\{a_\alpha\}_{\alpha\in I})+\lim_\omega(\{b_\alpha\}_{\alpha\in I}).$

\begin{proposition}\label{pseudo-rank}
Define $N:E\to [0,1]$ the following way.
If $r_\alpha\in \en_D(A_\alpha)$,
$$N(\{r_\alpha\}_{\alpha\in I})=\lim_\omega\frac
{\dim_D(\ra(r_\alpha))}{\dim_D(A_\alpha)},$$
where $\dim_D$ denotes the dimension function on finite dimensional right $D$-
modules.
Then $N$ is a pseudo-rank function.
\end{proposition}
\proof
Clearly, $N(1)=1$ and $N(0)=0$.
 Let $\underline{r}=\{r_\alpha\}_{\ai},\un{s}=
\{s_\alpha\}_{\ai}.$

\noindent
Then $\dim_D(\ra(r\a s\a))\leq \dim_D(\ra(r\a))$
and $\dim_D(\ra(r\a s\a))\leq \dim_D(\ra(s\a))$. 

\noindent
Hence $N(\un{r}\,\un{s})\leq N(\un{r}), N(\un{r}\,\un{s})\leq N(\un{s}).$
Now let $\un{e}=\{e\a\}_{\ai}, \un{f}=\{f\a\}_{\ai}\in E$ be orthogonal
idempotents.
Then for any $\alpha\in I, e\a$ and $f\a$ are orthogonal idempotents in
$\en_D(A\a)$. Thus
$$\dim_D(\ra(e\a+f\a))=\dim_D(\ra(e\a))+\dim_D(\ra(f\a))\,.$$
Consequently, $N(\un{e}+\un{f})=N(\un{e})+N(\un{f})\,.$\qed.

\noindent
Note that $\ke(N)=\{\un{x}\in E\,|\, N(\un{x})=0\}$ is a two-sided ideal of
$E$.
\begin{proposition} \label{continuous-regular-ring-construction} 
The ring $R_N=E/\ke (N)$ is a simple continuous von Neumann regular ring.
\end{proposition}
\proof
The direct product ring $E$ is a right and left self-injective von Neumann
regular ring, therefore by Corollary 13.5 \cite{Good} $E$ is continuous. If
$M$ is a maximal two-sided ideal of $E$, then by Corollary 13.27 \cite{Good}
$E/M$ is a simple continuous von Neumann regular ring. Thus it is enough to
show that $\ke (N)$ is a maximal ideal.
Let $Q\subset E$ be the following two-sided ideal: $\{a\a\}_{\ai}\in Q$
if $\{\alpha\in I\,|\,q_\alpha=0\}\in\omega\,.$ Then of course $E/Q$ is just
the ultraproduct of the rings $\en_D(A\a)$. Note that $Q$ is a prime ideal since
the ultraproduct of prime rings is a prime ring as well.
By Corollary 9.15 \cite{Good}, $E$ satisfies the 
general comparability axiom. 
Hence by Corollary 8.21 \cite{Good} $Q\subset \ke(N)$ is
a prime ideal. Proposition 16.25 \cite{Good} immediately implies that $\ke(N)$
is in fact a maximal two-sided ideal of $E$. This completes the proof of the
proposition.
\qed

\section{Sofic groups}\label{secsofic}
\def\idx {{F,\epsilon}}
\def\uRRR {E}
\def\RRR  {R(G)}
\def\urank {N}
\def\toRRR {T}
In this section, we recall the notion of a sofic group
from \cite{ES}
and prove Theorem~\ref{main-theorem}.

\begin{definition}[From \cite{ES}.]
For a finite set $V$ let $\map(V)$ denote the monoid of self-maps of
$V$ acting on the left, the monoid operation is the composition of
self-maps.
We say that two elements $e,f\in\map(V)$ are {\em $\epsilon$-similar}
for a real number $\epsilon\in (0,1)$, 
if the number of points $v\in V$ with $e(v)\neq f(v)$ is at most
$\epsilon|V|$. We say that $e,f$ are $(1-\epsilon)$-different, 
if the number of points $v\in V$ with $e(v)= f(v)$ is less than
$\epsilon|V|$, i.e. if they are not $(1-\epsilon)$-similar. 
\end{definition}

\begin{definition}[From \cite{ES}.]
\label{sofic-group}
The group $G$ is {\em sofic} if for each number
$\epsilon\in(0,1)$ 
and any finite subset $F\subseteq G$ there exists a finite set $V$ and
a function $\phi:G\to \map(V)$ with the following properties:
\begin{enumerate}
\item[(a)]   For any two elements $e,f\in F$ the map $\phi(ef)$ is
             $\epsilon$-similar to $\phi(e)\phi(f)$.
\item[(b)]   $\phi(1)$ is $\epsilon$-similar to the identity map of $V$.
\item[(c)]   For each $e\in F\setminus\{1\}$ the map $\phi(e)$ is
	     $(1-\epsilon)$-different from the identity map of $V$. 
\end{enumerate}
\end{definition}

{\bf Remark.}
The origin of this notion is \cite{Gro}, where Gromov introduced the
concept of initially subamenable graphs.
The term ``sofic group'' is introduced in \cite{Wei}: a finitely
generated group is called sofic if its Cayley graph is initially
subamenable.
The above definition is taken from \cite{ES} (with right action is
replaced by left action, which is more appropriate here).
This is the formulation which suits our need the best. 
In the case of finitely generated groups the two definitions 
in \cite{Wei} and \cite{ES} are equivalent. 
Although we shall not use it later, for the sake of completeness we
prove their equivalence. 
To begin with, we recall some notation, 
and the definition from \cite{Wei}.

\begin{definition}[From \cite{Wei}.]
\label{another-sofic}
Let $G$ be a finitely generated group, and $B\subset G$ a fixed
finite, symmetric (i.e. $B=B^{-1}$) generating set.
The Cayley graph of $G$ is a directed graph $\Gamma$, 
whose edges are labeled by the elements of $B$:
the set of vertices is just $G$, and the edges with label $b\in B$ are
the pairs $(g, bg)$ for all $g\in G$. 
Let $N_r$ denote the $r$-ball around $1\in\Gamma$ 
(it is an edge-colored graph, and also a finite subset in $G$). 
The group $G$ is called sofic in \cite{Wei}, if for
each $\delta>0$ and each $r\in\bN$ there is a finite directed graph
$(V,E)$ edge-labeled by $B$, and a subset $V_0\subset V$ with the
properties, that:
\begin{enumerate}
\item For each point $v\in V_0$ there is a function
	$\psi_v:N_r\to V$ which is an isomorphism (of labeled graphs)
	between $N_r$ and the $r$-ball in $V$ around $v$.
\item $|V_0|\ge(1-\delta)|V|$.
\end{enumerate}
\end{definition}

\begin{proposition}
For a finitely generated group $G$ the above two notion of soficity is
equivalent. In particular, Definition~\ref{another-sofic} does not
depend on the choice of the generating set $B$.
\end{proposition}

{\bf Proof of \ref{another-sofic}$\Rightarrow$\ref{sofic-group}:} 
Let $\epsilon>0$ and $F\subseteq G$ a finite subset. 
We chose $r\in\bN$ such that the product set $F\cdot F$ is contained
in $N_r$.
Let $(V,E)$ and $V_0\subset V$ be the labeled directed graph, and
subset corresponding to $\delta=\epsilon$ and $r$.
We shall use this finite set $V$, and 
define the function $\phi:G\to\map(V)$ as follows. 
For $g\in N_r$ and $v\in V$ let $\phi(g)(v)=\psi_v(g)$.
Otherwise, for $g\in G\setminus N_r$, 
we define $\phi(g)$ arbitrarily. It is an easy calculation to check
conditions (a)-(c) of Definition~\ref{sofic-group}.
\qed
\vskip 17pt
{\bf Proof of \ref{sofic-group}$\Rightarrow$\ref{another-sofic}:}
Let $\delta>0$ and $r\in\bN$. 
We set $F=N_{2r+2}$, and choose any $\epsilon>0$.
Let $\phi:G\to\map(V)$ be the function of
Definition~\ref{sofic-group} for this $(F,\epsilon)$.
We use this $V$ as the vertex set of our new graph, 
and for each $v\in V$ we define 
$\psi_v:N_{r+1}\to V$, $\psi_v(g)=\phi(g)(v)$.
Let $V_0$ be the set of those $v\in V$ for which 
\begin{enumerate}
\item[(A)]	$\psi_v(bg) =\psi_{\psi_v(g)}(b)$ 
		for all $g\in N_r$ and all $b\in B$,
\item[(C)]	$\psi_v(g)\ne\psi_v(h)$ whenever $g,h\in N_{r+1}$.
\end{enumerate}
Finally we build the labeled edges of $V$:
for each $b\in B$ and $v\in V$ we add the edge
$(v,\psi_v(b))$ with label $b$.
It is easy to see, that the $r$-ball around $v$ is the set
$\psi_v(N_r)$, and all edges coming out of it's points 
are contained in the $r+1$-ball $\psi_v(N_{r+1})$.
Condition (C) implies that $\psi_v$ is injective, and condition (A)
ensures that $\psi_v$ preserves the edges coming out of $\psi_v(g)$.
Hence Definition~\ref{another-sofic}~(1) is satisfied.
There are $|N_r|\cdot|B|$ equations to check in (A),
and there are $|N_{r+1}|^2$ inequalities in (C). 
We know from Definition~\ref{sofic-group}~(a) that each of the
equations can fail on at most $\epsilon|V|$ exceptional $v$.
Moreover, we shall apply Definition~\ref{sofic-group}
to the elements $g,h$ in (C), and we get
$$
\phi(g^{-1})(\psi_v(g))      \stackrel{\rm(a)}{=} 
\phi(g^{-1}g)(v)             \stackrel{\rm(b)}{=} 
v                            \stackrel{\rm(c)}{\ne}
\phi(g^{-1}h)(v)             \stackrel{\rm(a)}{=}
\phi(g^{-1})(\psi_v(h))
$$
(we marked each $=$ and $\ne$ with the corresponding condition of
Definition~\ref{sofic-group}), hence $\psi_v(g)\ne\psi_v(h)$
for all $v$ outside a subset of size at most $4\epsilon|V|$.
Hence each inequality in (C) can fail on at most $4\epsilon|V|$ 
exceptional $v$.
Hence Definition~\ref{another-sofic}~(2) holds if we choose
$\epsilon<\frac{\delta}{4|N_{r+1}|^2+|N_r|\cdot|B|}$.
\qed

For the sake of completeness, we quote without proof some important
properties of sofic groups:
\begin{proposition}[From \cite{ES}.]
\label{properties-of-sofic-groups}
Direct product, subgroup, inverse limit, direct limit, and free
product of sofic groups is sofic. 
If $N\lhd G$, $N$ is sofic and $G/N$ is amenable, then $G$ is also sofic. 
\end{proposition}

\begin{proposition}[From \cite{Wei} and \cite{ES}.]
\label{residually-amenable-is-sofic}
If $G$ is a locally residually amenable group then $G$ is sofic.
In particular, amenable and residually finite groups are sofic.
\end{proposition}
\proof
This proposition is proved in \cite{Wei} for finitely generated $G$,
and the general case follows from the fact, that soficness is a
property of finite subsets of $G$. 
It is also proved in \cite{ES}, it is an easy consequence of the above
Proposition~\ref{properties-of-sofic-groups}. 
For the convenience of the reader we sketch a proof here.

So let $G$ be a locally residually amenable group, and let $F\subset
G$ be a finite subset. Let $H$ be the subgroup generated by $F$, and
let $\bar H$ be an amenable factor group of $H$ such that the factor
map $F\to \bar F\subset H$ is a bijection. For elements $h\in H$ let
$\bar h\in\bar H$ denote the image of $h$. 
By Folner's theorem
there exists a finite subset $V\subseteq\bar H$ containing $\bar F$,
whose $\bar F$-boundary $\{v\in V\ |\ \bar Fv\nsubseteq V\}$ has at
most $\epsilon|V|$ elements. Next we define the function
$\phi:G\to\map(V)$:
$$
\phi(g)(v) = \left\{
{\bar gv \atop v}\ \ \ 
{\mbox{if }g\in H,\ \bar gv\in V \atop \mbox{otherwise}\hfill}
\right. 
$$
It is an easy calculation (left to the reader) to check
(a)-(c) of Definition~\ref{sofic-group}.
\qed

\vskip 17pt\noindent{\bf Proof of Theorem~\ref{main-theorem}:}
Thus let $G$ be a sofic group.
We define the index set:
$$
I = \left\{ (\idx) \ \Big|\ 
	F\subseteq G \mbox{ finite, and }
	\epsilon\in(0,1) \right\}
$$
and for each index $(H,\delta)\in I$ we define the nonempty subset
$$
I_{H,\delta} = \left\{ (\idx)\in I \ \Big|\ 
	H\subseteq F \mbox{ and }\epsilon\le\delta \right\}
\subseteq I
$$
The collection of nonempty subsets $\{ I_{H,\delta}\}$ is closed
under intersection, so there is an ultrafilter $\omega$ of subsets of
$I$ containing all $I_{H,\delta}$. Next, for each index
$(\idx)\in I$ we choose a finite subset $V_\idx$ and a
function $\phi_\idx:G\to\map(V_\idx)$ satisfying the
conditions (a)-(c) of Definition~\ref{sofic-group}. 
As in Section~\ref{omega-stuff}, 
for each index $\alpha\in I$ let $A_\alpha$ denote the right
$D$-module with basis $V_\alpha$,  $\uRRR$ will denote the
$\omega$-ultra-product of the endomorphism rings $\en_D(A_\alpha)$, 
$\urank$ is the pseudo-rank of Proposition~\ref{pseudo-rank},
and finally $\RRR=R_N$ will denote  factor ring constructed in
Proposition~\ref{continuous-regular-ring-construction}. 
This $\RRR$ is a simple, continuous von Neumann regular ring, 
this will be the ring we seek in the Theorem.
The elements of $\map(V_\alpha)$ extend linearly to
endomorphisms of the module $A_\alpha$, 
hence the functions $\phi_\alpha$ induce functions
$G\to\en_D(A_\alpha)$, and can be extended to linear functions
$$
T_\alpha:D(G)\to\en_D(A_\alpha)
$$
Taking the  ultra-product of these $T_\alpha$, and then composing with
the factor map $\uRRR\to\RRR$, we get a linear function:
$$
\toRRR :D(G)\to \RRR
$$
We shall prove that this $\toRRR$ is an injective homomorphism,
this will complete the proof of Theorem~\ref{main-theorem}. 
We see from property (a) of Definition~\ref{sofic-group} that for each $g,h\in G$
\begin{eqnarray*}
\urank\Big(\toRRR(g)\toRRR(h)-\toRRR(gh)\Big) &=&
\lim_\omega\frac{\r\Big(
	T_\alpha(g)T_\alpha(h)-T_\alpha(gh)\Big)}
   {\dim_D(A_\alpha)} \\
&\le&
\lim_\omega\frac{\left|\left\{v\in V_\alpha\Big|
	\phi_\alpha(g)(\phi_\alpha(h)(v))\ne \phi_\alpha(gh)(v)
	\right\}\right|}
   {|V_\alpha|} \\
&\le&
\lim_{\epsilon\to0} \epsilon = 0
\end{eqnarray*}
Therefore the map $\toRRR$ is a ring homomorphism.
The only that thing remains to be shown is the injectivity of
$\toRRR$.
So let $S\subset G$ be a finite subset, and for each $s\in S$ let
$k_{s}\in D$ be a nonzero element, we shall show that 
$\toRRR\Big(\sum_{s\in S} k_ss\Big)=
\sum_{s\in S} k_s\toRRR(s)\neq 0$ in $\RRR$.
For each index $\alpha=(\idx)\in I$ we choose a maximal subset 
$X_\alpha\subset V_\alpha$ such that
if $p,q\in X$ and $s,t\in S$ with either $p\ne q$ or $s\ne t$ then
$T_\alpha(s)(p)\ne T_\alpha(t)(q)$.
Since $X_\alpha$ is maximal, for each $p\in V_\alpha$ 
there is an element $q\in X_\alpha$ and elements $s,t\in S$
such that
$$
T_\alpha(p)(s)= T_\alpha(q)(t)
$$
The right hand side of this equation can take at most
$|S|\cdot|X_\alpha|$ different values. On the other hand,
by Definition~\ref{sofic-group}~(b) the map $T_\alpha(s)$ is injective on a subset
of size $(1-\epsilon)|V_\alpha|$, hence for each value of $s$ there are
at most $\epsilon|V_\alpha|+|S|\cdot|X_\alpha|$ possible value of
$p$. Since $p$ was arbitrary, we get
$|V_\alpha|\le\epsilon|S|\cdot|V_\alpha|+|S|^2\cdot|X_\alpha|$, or
$$
|X_\alpha| \ge \frac{1-\epsilon|S|}{|S|^2}|V_\alpha|
$$
If a non-zero element $x$ of $A_\alpha$ is spanned by $X_\alpha$ then
$\left(\sum_{s\in S} k_s T_\alpha(s)\right) (x)\neq 0\,.$
Hence
$$
\urank\left(\sum_{s\in S} k_s\toRRR(s)\right) \ge
\lim_\omega\frac{|X_\alpha|}{|V_\alpha|} \ge
\lim_{\epsilon\to0} \frac{1-\epsilon|S|}{|S|^2} = 
\frac{1}{|S|^2} >0
$$
This proves that $\toRRR\Big(\sum_{s\in S} k_ss\Big)\ne 0$.
Theorem~\ref{main-theorem} is proved now.
\qed
\begin{corollary}
If $G$ is a sofic group and $D$ is a division ring then $D(G)$ is
stably finite.
\end{corollary}

\end{document}